\documentclass[reqno,12pt]{amsart}
\usepackage{amssymb,epic,eepic,epsfig,graphics,graphicx,wrapfig,latexsym,color,enumerate,amsmath,amsthm,amsfonts}
\usepackage{a4wide}

\newtheorem{theorem}{Theorem}
\newtheorem{observation}{Observation}[section]
\newtheorem{proposition}{Proposition}[section]
\newtheorem{definition}[proposition]{Definition}
\newtheorem{notation}[proposition]{Notation}

\newtheorem{remark}[proposition]{Remark}
\newtheorem{lemma}[proposition]{Lemma}
\newtheorem{corollary}[proposition]{Corollary}

\newcommand{\stacksign}[2]{{\stackrel{\textrm{\scriptsize #1}}{#2}}}

\def\Tr{{{\rm Tr}}}
\def\rot{\mathrm{rot}}

\begin{document}

\title{The enumeration of planar graphs via Wick's theorem}

\author{Mihyun Kang$^{1}$}
\address{$^{1}$Humboldt-Universit\"{a}t zu Berlin, Institut f\"{u}r Informatik\\Unter den Linden 6, 10099 Berlin, Germany}
\email{kang@informatik.hu-berlin.de}

\author{Martin Loebl$^{2}$}
\address{$^{2}$Dept.~of Applied Mathematics and\\
Institute of Theoretical Computer Science (ITI)\\
Charles University \\
Malostranske n. 25 \\
118 00 Praha 1 \\
Czech Republic
\newline
{\em and}  Depto. Ing. Matem\,atica \\ University of Chile \\ Chile.
\newline
M.L.gratefully acknowledges the support of CONICYT via Anillo en Redes ACT08}
\email{loebl@kam.mff.cuni.cz}

\maketitle

\begin{abstract}
A seminal technique of theoretical physics called {\em Wick's theorem} interprets the Gaussian
matrix integral of the products of the trace of powers of Hermitian matrices 
as the number of labelled {\em maps} with a given degree sequence, sorted by their
Euler characteristics. This leads to the map enumeration results analogous to those obtained
by combinatorial methods. In this paper we show that the enumeration
of the graphs {\em embeddable} on a given $2$-dimensional surface (a main research topic of contemporary enumerative combinatorics) can also be formulated as the Gaussian matrix integral of an ice-type partition function. 
Some of the most puzzling conjectures of discrete mathematics are related to the notion of the {\em cycle double cover}. 
We express the number of the graphs with a fixed {\em directed cycle double cover}  as the Gaussian matrix integral of an Ihara-Selberg-type function.

\end{abstract}

\section{Introduction}
\label{sec.int}
Let $M=(M_{ij})$ be an $N\times N$ Hermitian matrix, i.e., $M_{ij}={\overline M}_{ji}$ for every $1\leq i,j\leq N$, where ${\overline M}_{ji}$ denotes the complex conjugate of ${M}_{ji}$. Let $\varphi(M)= \sum_I a_I\prod_{ij\in I}M_{ij}$ be a polynomial in its entries, where $I$ ranges over a finite system of multisets of elements
of $N\times N$ and $a_I$'s are complex numbers. We start with the following notation. 

\begin{notation}
\label{def.uno}
Denote $\langle \cdot \rangle $ by
$$
\langle \varphi \rangle \; =\;  \sum_I a_I \langle \prod_{ij\in I}M_{ij} \rangle 
\; =\;  \sum_I a_I\sum_P \prod_{(p,q)\in P} \langle M_pM_q \rangle ,
$$
where $P$ ranges over all partitions of $I$ into pairs, and for $p=(p_1,p_2), q=(q_1,q_2)$
we have $\langle M_pM_q \rangle $ is non-zero if and only if $p_1=q_2$ and $p_2=q_1$ and in that case $\langle M_pM_q \rangle \; =\;  1/N$.
\end{notation}

\bigskip
\subsection{The motivation}
\label{sub.mot}
Wick's theorem interprets $\langle \varphi \rangle $ as the Gaussian matrix integral (see Section \ref{sec.wick}).  
This has been successfully used in enumeration of combinatorial structures and in other fields of mathematics~\cite{BIZ,BFG,BFG03,BIPZ,DF,DF04,HZ86,Kontsevich92,LandoZvonkin,MulaseWaldron,Okounkov,Penner,PS06,Zvonkin}.

The basic combinatorial structure studied in this paper is that of a {\em graph}. A graph $G$ is a pair $(V,E)$ where $V$ is the set of {\em vertices} and $E$ is the set of pairs of vertices called {\em edges}. The vertices
of an edge are called its {\em endvertices}. A graph may have {\em multiple edges}, i.e., different edges with 
the same endvertices, and also {\em loops}, i.e., edges whose two endvertices are equal. The graphs without 
multiple edges and loops are called {\em simple}. A {\em map} is a graph together with a fixed cyclic ordering of the incident edges of each vertex: it defines an embedding of the graph on an orientable $2$-dimensional surface (see~\cite{mohartomassen}). Note that a map is a graph that is {\it embedded} on a given $2$-dimensional surface. 

Principally studied functions $\varphi$ are products of powers of the trace of $M^i$ for a nonnegative integer $i$. For such functions $\varphi$, $\langle \varphi\rangle $ has a useful graphic interpretation as the number of labelled {\it maps} with given degree sequence, sorted by their Euler characteristics (these maps are the {\it Feynman diagrams} for the matrix integral). This 
seminal formula was first described in \cite{BIPZ}. 
While the Gaussian matrix integral of a polynomial in the traces of the powers of $M$ is well understood, the Feynman diagrams expansions of the Gaussian matrix integral of the {\em exponential} of such polynomials hold in general only on the level of {\em formal power series}. We will denote this formal equality by $\ '=\, '$.
The expression $f(x) \ '=\, ' \ g(x)$ simply means that all the derivatives of $f,g$ are equal when the variable $x$ is set to zero. Similarly the expression $f(x)\ '\leq\, '\  g(x)$ means that all the derivatives of $f,g$ satistfy the inequality when the variable $x$ is set to zero. Finally, let $s(N)$ be a positive integer function in $N$ and let $f_N(x), g_N(x)$ be two formal power series in variable $x$, whose coefficients depend on a positive integer $N$. 
The expression $f_N(x)\ '=_{s(N)}\, '\  g_N(x)$ (resp. $f_N(x)\ '\leq_{s(N)}\, '\  g_N(x)$) means that all the $k$-th derivatives of $f,g$, for $k\leq s(N)$, satistfy the equality (resp. inequality) when the variable $x$ is set to zero.
These notations are defined analogously for the multivariate generating functions. We describe the classical application of the matrix integrals in the next subsection.

\bigskip
\subsection{Enumeration of maps}
\label{sub.mapps}
Let $\mathcal H_N$ denote the set of all $N\times N$ Hermitian matrices. 

We consider a function $\psi$ which maps $M\in \mathcal H_N$ to
\begin{eqnarray*}
\psi(M):=e^{(-N)\sum_{i\in \mathbb N}z_i\Tr\left(\frac{M^i}{i}\right)}, \label{func.map}
\end{eqnarray*}
where $\mathbb N$ is the set of all positive integers and $z_i$'s are formal variables.

Taking the logarithm of $\langle \psi\rangle $ we get the following formula for connected maps sorted by genus as discussed by Br\'ezin, Itzykson, Parisi and Zuber in \cite{BIPZ}.

\begin{theorem}[\cite{BIPZ}]\label{thm.maps}
\begin{eqnarray*}
&&\log \langle \psi\rangle  \; \ '=\, '\  \; \sum_{g\ge 0}N^{2-2g} \sum_{(n_1,\cdots,n_k)\in \mathbb N^k}\prod_{i=1}^k\frac{(-z_i)^{n_i}}{n_i!}M_g(n_1,\cdots,n_k)
\end{eqnarray*}
where $M_g(n_1,\cdots,n_k)$ denotes the number of connected maps with genus $g$ and $n_i$ vertices of degree $i$ for $1\le i\le k$.
\end{theorem}

The convergence of the generating function in Theorem \ref{thm.maps} was studied extensively~\cite{EM,EMP,Gui,GuiICM,GuiMau,Mulase}. For example, the following theorem was proved by Ercolani and McLaughlin in \cite{EM}. (For more general results, see a lecture note by Guionnet~\cite{GuiLecture} and references therein.) 

\begin{theorem}[\cite{EM}]
\label{thm.emgui}
Let a function $\tilde \psi$ map each Hermitian matrix $M\in \mathcal H_N$ to 
$$\tilde \psi(M):=e^{(-N)\sum_{i\in I}z_i\Tr\left(\frac{M^i}{i}\right)},$$
where $I$ is a finite set of positive integers such that its maximum element is even.
Then
for each $i\in I$ there is $\epsilon_i>0$ so that for $z_i$ a non-zero real variable with $z_i\in (0,\epsilon_i)$,
Theorem \ref{thm.maps} holds as equalities between {\em analytic functions} (in variables $z_i$) when
$\psi$ is replaced by $\tilde \psi$.
\end{theorem}

\bigskip
\subsection{Main contribution}
\label{sub.mainn}

A natural question is whether the Gaussian matrix integral method may as well be applied to the enumeration of graphs that are {\it embeddable} on a given $2$-dimensional surface. This enumeration has been a long-standing open problem in the enumerative combinatorics, solved only for planar graphs and fairly recently by Gim{\'e}nez and Noy~\cite{GN}.

One of our main results (Theorem \ref{thm.main1}) is that this enumeration task 
can as well be formulated as a Gaussian matrix integral of an ice-type partition function 
(but whose asymptotic analysis is not carried out so far).

Another main result (Theorem \ref{thm.main2}) is that the number of graphs with a fixed {\it directed cycle double cover} is expressed as a Gaussian matrix integral of an Ihara-Selberg-type function. 
Even though the study of the cycle double cover conjectures is extensive, no attempt has been made before to calculate these numbers.

The main contribution of this paper is in applying the matrix integral
to {\it combinatorially defined functions}, in order to loosen the strong connection with maps implied by 
the integration of the functions of traces.
In the rest of this section we will show how to achieve this goal, introduce necessary concepts, and state our main results more precisely.

\bigskip
\subsection{Directed graphs}
\label{sec.di}
A starting idea is to associate  directed graphs with matrices. 
Let $M\in \mathcal H_N$ and let $D= D(M)= (N, N\times N)$ be a complete directed graph with weights on directed edges given by $M$.
\begin{definition}
\label{def.closedwalk}
A {\em closed walk} in $D$ is defined as a circular sequence of directed edges of $D$ so that the tail of each edge equals the head of the consecutive edge. A closed walk is called {\em pointed} if it has a prescribed beginning. A {\em closed trail} is a closed walk which contains each directed edge {\em at most once}.
\end{definition}
Now $\Tr(M^n)$ can be interpreted as
$\sum_p\prod_{e\in p}M_e$, where the sum is over all pointed closed walks in $D$ of length $n$.
Similarly, $\Tr(M^3)^4\Tr(M^2)^3$ can be interpreted as  a function which maps each Hermitian matrix $M\in \mathcal H_N$ to  
$$
(\sum_{p_1}\prod_{e\in p_1}M_e)^4(\sum_{p_2}\prod_{e\in p_2}M_e)^3,
$$
where the first sum is over all pointed closed walks $p_1$ in $D$ of length $3$, and the second sum is over all pointed closed walks $p_2$ in $D$ of length $2$.

\begin{definition}
\label{def.eul}
 A subset $A$ of directed edges of $D$ is called {\em eulerian}
if for each vertex its indegree in $A$ is equal to its outdegree in $A$. 
We denote by ${\mathcal E}=\mathcal E(D)$ the collection of all eulerian sets in $D$. 
\end{definition}

\begin{definition}
\label{def.even}
A subset $A$ of directed edges of $D$ is called {\em even} if $A$  can be written as a union of edge-disjoint directed cycles  of length bigger than two. 
\end{definition}

For a set $A$ of directed edges of $D$,  a subset $C\subset A$ is called a {\em component} of $A$
if $C$ is a maximal subset of $A$ with respect to  inclusion that induces a connected underlying undirected graph.

\bigskip
\subsection{Enumeration of embeddable graphs}
\label{sec.embb}
We say that an integer function $f(N)$ is of type $O(N^{-1})$ if $f(N)$ is bounded above by $cN^{-1}$ for a constant $c$ and all large enough $N$.

\smallskip
Let $\eta$ be a function which
maps each Hermain matrix $M\in \mathcal H_N$ to
\begin{eqnarray}
\eta(M,Nzy^{-1},y)= \sum_{A\in \mathcal E}(\prod_{e\in A}M_e)\prod_{C:\text{ component of } A}\frac{(Nz)^{|C|/2+1}-Nzy^{|C|/2}}{Nz-y}.\label{def.eta}
\end{eqnarray}

\bigskip
The first result of the paper is the following expression for planar graphs.
\begin{theorem}
\label{thm.main1}
Let $s(N)\le \frac{N}{2}$ be any function in $N$ which approaches to $\infty$ slower than $N$ when $N\to \infty$.
The Gaussian integral of the function $\eta$ defined in (\ref{def.eta}) satisfies the following.
\begin{eqnarray*}
&&\exp\left(\sum_{n\ge 1} \sum_{r\ge 0 }\left[\, p(n,r)+e_1(N,n,r)\, \right]\frac{z^{r}y^{n-2}}{n!}\right)\ '\leq_{s(N)}\, ' \\
&&e^{-N^2}\langle  \eta(M,Nzy^{-1},y) \rangle \ \\
&&'\leq_{\frac{N}{2}}\, ' \ \exp\left(\sum_{n\ge 1} \sum_{r\ge 0 }\left[\, p(n,r)+e_2(N,n,r)\, \right]\frac{z^{r}y^{n-2}}{n!}\right),
\end{eqnarray*}
where $p(n,r)$ denotes the number of labelled connected simple graphs on $n$ vertices which have {\em planar} embeddings with $r$ faces, and $e_1(N,n,r), e_2(N,n,r)$ are functions of type $O(N^{-1})$.
\end{theorem}

Let $p(n)$ denote the number of labelled planar connected simple graphs on $n$ vertices. Clearly, we have that $p(n)= \sum_{r\ge 0 }p(n,r)$ and that  ${p(n)}/{n!}$ is equal to the coefficient of $y^{n-2}$ in 
$\lim_{N\to \infty}\frac{1}{N^2}\log\langle  \eta(M,Nzy^{-1},y) \rangle |_{z=1}$.

It would be very interesting to see if we can recover 
$$p(n)\sim c\ n^{-7/2}\ \gamma^n \ n!$$ with $c\doteq 0.478$ and $\gamma\doteq 27.2$, which is computed by Gim\'emez and Noy~\cite{GN}.

The embeddable graphs of higher genus contribute to the subleading terms \linebreak $e_1(N,n,r), e_2(N,n,r)$ in  Theorem \ref{thm.main1}; it should be possible to derive their enumeration from our formulas in Section \ref{sec.planargraph}.

The function $\eta$ in (\ref{def.eta}) is an ice-type partition function. Ice-type partition functions play a seminal role in modern
mathematics and theoretical physics (see \cite{Wu}), and we hope that the asymptotic analysis of the
matrix integral of the function $\eta$ is not out of reach.

\bigskip
\subsection{Enumeration of directed cycle double covers}
\label{sec.dcdc}
Observe that an eulerian set does not need to be even.
Even sets are closely related to the {\em cycle double cover conjectures}. 

\begin{definition}
\label{def.DDCC}
Let $G$ be an undirected graph. A collection of its cycles is called a {\em cycle double cover (CDC)}
if each edge belongs to exactly two of the cycles. Moreover it is called a {\em directed cycle double cover
 (DCDC)} if it is possible to orient the cycles so that they go oppositely on each edge.
\end{definition}

Some of the most puzzling conjectures of discrete mathematics are centred around this notion. A graph is {\em bridgeless} if it cannot be disconnected by deletion
of a single edge. Clearly a graph with a bridge does not have a CDC.
On the other hand, there are 
\begin{itemize}
\item
{\em Cycle double cover conjecture}: Is it true that each bridgeless graph has a CDC?
\item
{\em Directed cycle double cover conjecture}: Is it true that each bridgeless graph has a DCDC?
\end{itemize}

Our next main result concerns directed cycle double covers.

\begin{theorem}
\label{thm.main2}
Let ${\mathcal A}$ be the set of all pairs $(q,K)$ where $q$ is an even set of edges of $D$
and $K$ is a decomposition of $q$ into directed cycles of length at least three.
Then
$$
\langle \sum_{(q,K)\in {\mathcal A}}\prod_{e\in q}M_e\rangle \; = \; \sum_{[(G,C)]}\frac{N(N-1)\dots (N-|V(G)|+1)}{|{\rm Aut}(G,C)|N^{e(G)}},
$$
where the sum is over all isomorphism classes of pairs $(G,{C})$ where $G$ is a simple graph with at most $N$ vertices and no isolated vertex, and  ${C}$ is a {\em specified} DCDC of $G$.
\end{theorem}

We show in Section \ref{subsec.ihara} that the integrand $\xi$ in Theorem \ref{thm.main2} which maps $M\in \mathcal H_N$ 
to 
\begin{eqnarray}
\xi(M)\; := \; \sum_{(q,K)\in {\mathcal A}}\prod_{e\in q}M_e \label{fun.dcdc}
\end{eqnarray}
is an Ihara-Selberg-type function (see \cite{ST1}, \cite{ST2}). Thus again,
we hope that an analysis of the matrix integral of the function $\xi$ may be achievable.

\bigskip

The rest part of the paper is organized as follows. In Section \ref{sec.wick} we review the Gaussian matrix integral method and Wick's theorem. In Section \ref{sec.wdc} we introduce identities essential to the enumeration of embeddable graphs on a 2-dimensional surface. In Sections \ref{sec.planargraph} and \ref{sec.dcdc} we prove our main Theorems \ref{thm.main1} and \ref{thm.main2}.


\bigskip
\section{Gaussian matrix integral and Wick's theorem}
\label{sec.wick}
In this section, based on \cite{DF04}, we recall the definition of Gaussian matrix integral, useful identities, and Wick's theorem in particular. 

Let $M\in \mathcal H_N$ and let $$dM=\prod_i dM_{ii} \prod_{i<j} d\; {\rm Re}(M_{ij}) d\; {\rm Im}(M_{ij})$$ denote the standard Haar measure on $\mathcal H_N$, where ${\rm Re}(M_{ij})$ and ${\rm Im}(M_{ij})$ denote the real part and the imaginary part of $M_{ij}$.

For an arbitrary function $\varphi:\mathcal H_N\to \mathbb C$, the Gaussian matrix integral of  $\varphi$  is defined as 
\begin{eqnarray}
\langle \varphi\rangle \;& =&\; \frac{1}{Z_0(N)}\int_{\mathcal H_N} e^{-N\; \Tr(\frac{M^2}{2})}\varphi(M) dM,\label{def.matin}
\end{eqnarray}
where $Z_0(N)$ is the normalization factor making $\langle 1\rangle \; =1$, that is,
$$Z_0(N)=\int_{\mathcal H_N} e^{-N\; \Tr(\frac{M^2}{2})}dM.$$

We are particularly interested in a function $\varphi$ which maps $M\in \mathcal H_N$ to 
\begin{eqnarray}
 \varphi(M)=\sum_I a_I\prod_{ij\in I}M_{ij},\label{eq.varphi}
\end{eqnarray}
where $I$ ranges over a finite system of multisets of elements of $N\times N$. To compute $\langle \varphi \rangle$ we introduce  the {\em source integral} $\langle e^{\Tr(MS)} \rangle $ for a given $N\times N$ Hermitian matrix $S$, where $\Tr(MS)$ denotes the trace of the matrix $MS$.
It can easily be computed by 
\begin{eqnarray}
\langle e^{\Tr(MS)} \rangle 
\; &=&\; \frac{1}{Z_0(N)} \int_{\mathcal H_N} e^{-N\; \Tr(\frac{M^2}{2})} e^{\Tr(MS)} dM\nonumber\\
\; & =&\; \frac{1}{Z_0(N)} \int_{\mathcal H_N} e^{-N\; \Tr(\frac{1}{2}(M-\frac{S}{N})^2)}e^{\frac{\Tr(S^2)}{2N}} dM\nonumber\\
\; &=&\; e^{\frac{\Tr(S^2)}{2N}},\label{eq.matsource}
\end{eqnarray}
since the trace is linear and $\Tr(MS)=\Tr(SM)$.
On the other hand, for any $1\le i,j\le N$, 
\begin{eqnarray*}
\frac{\partial}{\partial S_{ji}}e^{\Tr(MS)}\Big|_{S=0}
\; &=&\; \left(\frac{\partial}{\partial S_{ji}} \Tr(MS)\right) e^{\Tr(MS)}\Big|_{S=0}\\
\; &=&\; \left(\frac{\partial}{\partial S_{ji}} \sum_{m,n}M_{mn}S_{nm}\right) e^{\Tr(MS)}\Big|_{S=0}\\
\; &=&\; M_{ij}.
\end{eqnarray*}
Thus the derivatives of the source integral becomes
\begin{eqnarray}
&&\frac{\partial}{\partial S_{ji}}\frac{\partial}{\partial S_{lk}}\cdots \langle e^{\Tr(MS)} \rangle \Big|_{S=0}\nonumber\\
\; &=&\; \frac{1}{Z_0(N)} \int_{\mathcal H_N} e^{-N\; \Tr(\frac{M^2}{2})} \frac{\partial}{\partial S_{ji}}\frac{\partial}{\partial S_{lk}}\cdots e^{\Tr(MS)}\Big|_{S=0} dM\nonumber\\
\;&=&\; \frac{1}{Z_0(N)} \int_{\mathcal H_N} e^{-N\; \Tr(\frac{M^2}{2})} M_{ij}M_{kl}\cdots dM\nonumber\\
\;&=&\; \langle  M_{ij}M_{kl}\cdots \rangle ,\label{eq.matder}
\end{eqnarray}
where the first equality is due to the Leibniz integral rule.

Using (\ref{eq.matder}) and (\ref{eq.matsource}), we obtain
\begin{eqnarray}
\langle  M_{ij}M_{kl}\cdots \rangle 
\; & &\stacksign{(\ref{eq.matder})}{=}\; 
\frac{\partial}{\partial S_{ji}}\frac{\partial}{\partial S_{lk}}\cdots \langle e^{\Tr(MS)} \rangle \Big|_{S=0}\nonumber\\
\; & &\stacksign{(\ref{eq.matsource})}{=}\; \frac{\partial}{\partial S_{ji}}\frac{\partial}{\partial S_{lk}}\cdots e^{\frac{\Tr(S^2)}{2N}}\Big|_{S=0}\label{eq.matpair}
\end{eqnarray}
and in particular
\begin{eqnarray}
\langle  M_{ij}M_{kl} \rangle 
\; &=&\; \frac{\partial}{\partial S_{ji}}\frac{\partial}{\partial S_{lk}} e^{\frac{\Tr(S^2)}{2N}} \Big|_{S=0}\nonumber\\
\; &=&\; \frac{\partial}{\partial S_{ji}}\left(\frac{\partial}{\partial S_{lk}}\frac{\Tr(S^2)}{2N}\right)e^{\frac{\Tr(S^2)}{2N}}\Big|_{S=0}\nonumber\\
\; &=&\; \frac{\partial}{\partial S_{ji}}\left(\frac{\partial}{\partial S_{lk}}\frac{\sum_{m,n}S_{mn}S_{nm}}{2N}\right)e^{\frac{\Tr(S^2)}{2N}}\Big|_{S=0}\nonumber\\
\; &=&\; \frac{\partial}{\partial S_{ji}}\frac{S_{kl}}{N}e^{\frac{\Tr(S^2)}{2N}}\Big|_{S=0}\nonumber\\
\; &=&\; \frac{\delta_{il}\delta_{jk}}{N}.\label{eq.matnonzero}
\end{eqnarray}
Further, it is clear that the derivatives in~(\ref{eq.matpair}) and ~(\ref{eq.matnonzero}) must be taken in pairs (e.g., $S_{ji}$ and $S_{lk}$ with $l=i$ and $k=j$) to get a non-zero contribution. This yields the following result known as \emph{Wick's theorem}.
\begin{theorem}[Wick's theorem]\label{thm.wicktheorem}
Let $M\in \mathcal H_N$ and $I$ be a multiset of elements of $N\times N$. Then
\begin{eqnarray}
\langle \prod_{ij\in I}M_{ij} \rangle 
\; &=&\; \sum_{\text{pairing } P\subset I^2}\quad \prod_{(ij,kl)\in P} \langle  M_{ij}M_{kl} \rangle \nonumber\\
\; &=&\; \sum_{\text{pairing } P\subset I^2}\quad\prod_{(ij,kl)\in P}\frac{\delta_{il}\delta_{jk}}{N}.\label{eq.wicktheorem}
\end{eqnarray}
\end{theorem}

Due to the linearity of the integral and (\ref{eq.wicktheorem}),  for the function $\varphi$ defined in (\ref{eq.varphi}) we have
\begin{eqnarray*}
\langle \varphi \rangle \; =\;  \sum_I a_I \langle \prod_{ij\in I}M_{ij} \rangle \; =\;  \sum_I a_I\sum_P \prod_{(p,q)\in P} \langle M_pM_q \rangle ,
\end{eqnarray*}
where $P$ ranges over all partitions of $I$ into pairs, and for $p=(p_1,p_2), q=(q_1,q_2)$
we have $\langle M_pM_q \rangle $ is non-zero if and only if $p_1=q_2$ and $p_2=q_1$ and in that case $\langle M_pM_q \rangle \; =\;  1/N$.
This is how we arrive at the notation of $\langle \cdot \rangle $ in Notation~\ref{def.uno}.


\bigskip
\section{Trail double cover and pseudoembedding}
\label{sec.wdc}
In this section we derive several identities, which are essential for applying the Gaussian matrix integral method to the enumeration of embeddable graphs on a 2-dimensional surface. 
We recall that a {\em map} is a graph together with a fixed cyclic ordering of the incident edges
of each vertex. This ordering (of each vertex) defines an embedding of the graph on an orientable $2$-dimensional surface
(see~\cite{mohartomassen}).
A map is also called a {\em fat graph}; we prefere this term since it corresponds to a
helpful graphic representation of Wick's theorem. In a fat graph $F$  the vertices are made into discs (islands) and connected by fattened
edges (bridges) prescribed by the cyclic orders. This defines a two-dimensional orientable surface with boundary which we also denote by $F$. 
Each component of the boundary of $F$ will be called a {\em face} of $F$. Each face
is an embedded circle (see e.g.,~\cite{mohartomassen}).
We will denote by $G(F)$ the underlying graph of $F$.

We denote by $e(F), v(F), f(F), c(F)$, and $ g(F)$ the number of the edges, vertices, faces,
connected components, and genus of $F$. 
We recall that $2g(F)= 2c(F)+ e(F)- v(F)- f(F)$.

In the next sections we will count fat graphs and their relatives. 
To avoid confusion we assume that a fat graph has {\em labelled} vertices,
i.e., two fat graphs are equal if they are equal as sets. We speak about
{\em unlabelled} fat graphs if the equality is up to isomorphism.

Let $M\in \mathcal H_N$ and we recall that $D= D(M)= (N, N\times N)$ is a directed graph
with weights on directed edges given by $M$. A closed walk in $D$ is defined as a circular sequence of directed edges of $D$ so that the tail of each edge equals the head of the consecutive edge. As we have seen in Section \ref{sec.di}, $\Tr(M^3)^4\Tr(M^2)^3$ can be interpreted as  
$$
(\sum_{p_1}\prod_{e\in p_1}M_e)^4(\sum_{p_2}\prod_{e\in p_2}M_e)^3,
$$
where the first sum is over all pointed closed walks $p_1$ in $D$ of length $3$, and the second sum is over all pointed closed walks $p_2$ in $D$ of length $2$.

A {\em proper pairing} of a subset $A$ of directed edges of $D$ is a partition of $A$ into pairs $ij,ji\in N\times N$ 
of oppositely directed edges.
Hence if $\varphi(M)\; = \; \Tr(M^3)^4\ \Tr(M^2)^3$ we get
$$
\langle \varphi \rangle \; = \; \sum_{q=q_1q_2 \dots q_7}\sum_P\prod_{(e,e')\in P}{1/N},
$$
where the second sum is over all {\it proper pairings} $P$ of the directed edges
of the disjoint union $q=q_1q_2 \dots q_7$ of $7$ pointed closed walks, from which four have length $3$ and remaining three have length $2$. Two directed edges form a proper pairing if one is reversed the other. We also say that such a pair $(q,P)$  {\em contributes} to $\langle \varphi \rangle $.

Next we  recall that a closed trail is a closed walk which does not repeat edges and we denote by ${\mathcal W}(r)$ the set of all subsets of edges of $D$ that can be decomposed into $r$ edge-disjoint closed trails,
and we let ${\mathcal W}= \cup_{r\geq 0}{\mathcal W}(r)$. 

We define a function $\omega_r$ which maps each Hermitian matrix $M\in \mathcal H_N$ to 
\begin{eqnarray}
\omega_r(M)\; := \; \sum_{q\in {\mathcal W}(r)}y^{|q|/2}\prod_{e\in q}M_e,\label{def.combfunc}
\end{eqnarray}
for each integer $r\ge 1$ and $\omega_0(M)=1$.
We call $y^{|q|/2}\prod_{e\in q}M_e$ the weight of $q$.

Next two definitions and proposition are crucial.

\begin{definition}
\label{def.waco}
Let $G$ be a graph. A {\em trail double cover} (TDC) is a collection of closed trails so that each edge of $G$ is traversed
exactly twice, and in the opposite directions.
\end{definition}

We remark that each graph has a TDC.
Next we introduce a useful notation. A simple graph will be called {\em nimple} if it has no vertex of degree 0. 
\begin{definition}
\label{def.we}
Let $G$ be a finite nimple graph with at most $N$ vertices. 
Then let $c(G)$ be the set of all pairs $(q,P)$ so that there is a colouring $d$ of the vertices 
of $G$ by colours $\{1,\dots, N\}$, where each vertex gets a different colour, $$q= \{(d(x),d(y)),(d(y),d(x));\{x,y\}\in E(G)\},$$
and $P$ consists of all the pairs $\{[(d(x),d(y)),(d(y),d(x))];\{x,y\}\in E(G)\}$. 
\end{definition}

We remark that each such $q$ belongs to ${\mathcal W}$ since $G$ has a TDC and that 
$$|c(G)|\; = \;  \frac{N(N-1)\cdots (N-|V(G)|+1)}{|{\rm Aut}(G)|}.$$

\begin{proposition}
\label{prop.la}
Let $r\ge 0$. A pair $(q,P)$  with $P$ a proper pairing of $q$ contributes to $\langle \omega_r \rangle $ if and only if
 there is a nimple graph $G$ with a TDC consisting of $r$ closed trails such that $(q,P) \in c(G)$.
\end{proposition}
\begin{proof}
If $(q,P)\in c(G)$, then any TDC consisting of $r$ closed trails provides a partition of $q$ into its trails 
and hence $(q,P)$ contributes to $\langle \omega_r \rangle $. 
On the other hand, if $(q,P)$ contributes to $\langle \omega_r \rangle $, then letting $G$ be the graph with the vertices from $\{1,\dots, N\}$ 
and the edges given by $P$, we get that $G$ is nimple since $q$ consists of edge-disjoint closed trails. 
Furthermore, $G$ has a TDC consisting of $r$ closed trails and hence $(q,P)\in c(G)$.
  \end{proof}

 \begin{proposition}
\label{prop.lla}
If $c(G)\cap c(G')\neq \emptyset$, then $G$ is isomorphic to $G'$. Moreover, if $G$ is isomorphic to $G'$, then $c(G)= c(G')$.
\end{proposition}
\begin{proof}
If $(q,P)\in c(G)\cap c(G')$, then the construction of $q$ induces a function between the sets of vertices of $G$ and $G'$, and $P$ gives the edges of both $G, G'$. Hence they are isomorphic.
The second part is true since the definition of $c(G)$ does not depend on 'names' of the vertices. 
  \end{proof}

As a consequence we obtain the following.

\begin{theorem}
\label{thm.main7}
Let $\omega_r$ be the function defined in (\ref{def.combfunc}). Then
$$
\langle \omega_r \rangle \; = \; \sum_{[G]}y^{e(G)}\frac{N(N-1)\cdots (N-|V(G)|+1)}{|{\rm Aut}(G)|N^{e(G)}},
$$
where the sum is over all isomorphism classes $[G]$ of nimple graphs with at most $N$ vertices that have a TDC consisting of $r$ closed trails.
\end{theorem}

Next we need to extend the notion of the embedding. A {\em pseudosurface} is obtained from a surface
by identifying finitely many (not necessarily disjoint) pairs of vertices. The identified vertices are
called {\em singularities} of the pseudosurface. 
A {\em pseudoembedding} of a graph $G$ is defined to be an embedding on a pseudosurface $S$ such that each singularity
of $S$ has a vertex of $G$ embedded in it.

\begin{proposition}
\label{p.embed}
Let $G$ be a nimple graph.
There is a TDC of $G$ consisting of $r$ closed trails iff there is a pseudoembedding
of $G$ with $r$ faces.
\end{proposition}
\begin{proof}
The first implication: we construct, from each trail $w$ of TDC, a planar disc bounded
by a polygon with $|w|$ edges. Let the boundary be directed according to the direction of $w$.
Then we glue the discs together by identifying the oppositely oriented edges with the same endvertices. 
The result is a pseudosurface.

The other implication: given a pseudoembedding, orient the faces so that each edge appears twice,
and with the opposite orientation. 
The orientations of the faces clearly form a TDC.

\end{proof}

As a consequence of Theorem \ref{thm.main7} and Proposition \ref{p.embed} we have

\begin{corollary}
\label{cor.main2}
Let $\omega_r$ be the function defined in (\ref{def.combfunc}). Then
$$
\langle \omega_r \rangle \; = \; \sum_{[G]}y^{e(G)}\frac{N(N-1)\dots (N-|V(G)|+1)}{|{\rm Aut}(G)|N^{e(G)}},
$$
where the sum is over all isomorphism classes $[G]$ of nimple graphs with at most $N$ vertices that have pseudoembeddings
with $r$ faces.
\end{corollary}

If $W$ is a fat graph, then its {\em dual} is the abstract graph $W^*= (V^*,E^*)$ with vertex set $V^*$ and edge set $E^*$, where $V^*$ is the set of (the duals of) the faces of $W$  and $E^*$ is the set of (the duals of) the edges $e$ of $W$ such that in $W^*$, $e=e(fg)$ has endvertices $f,g\in V^*$, and in $W$, $e$ lies on the boundary of faces $f,g$ of $W$. Note that there is a natural bijection between the set
of the edges of $W^*$ and that of $W$. 

Next we consider the geometric duals of the pseudoembeddings.

\begin{definition}
\label{def.rrrl}
Given a fat graph $W$  and  a partition $Q$ of the set of its faces,  we denote by
$G_Q(W^*)$ the (abstract) graph obtained from $W^*$ by the contraction of the classes
of $Q$ into single vertices.

A  pair $(W,Q)$, where $W$ is a fat graph and $Q$ a partition of the set
of its faces, is called {\em $r$-relevant} if $W$ has $r$ vertices, $Q$ has at most $N$ classes, and $G_Q(W^*)$ is a nimple graph.

We say that $(W,Q)$ is equivalent to $(W',Q')$, denoted by $(W,Q)\sim (W',Q')$, if $G_Q(W^*)$ is isomorphic to $G_{Q'}(W'^*)$.
\end{definition}

We call a partition $Q$ of the face set of $W$ {\em trivial} and denoted it by  $\emptyset$ if each face of $W$ forms a single partition class of $Q$.

\smallskip
The following is a consequence of Corollary \ref{cor.main2}.

\begin{corollary}
\label{cor.main3}
Let $\omega_r$ be the function defined in (\ref{def.combfunc}). Then
\begin{eqnarray}
\langle \omega_r \rangle \; = \; \sum_{[(W,Q)]}y^{e(W)}\frac{N(N-1)\dots (N-\alpha(Q)+1)}{|{\rm Aut}(G_Q(W^*))|N^{e(W)}},\label{eq.equivalent}
\end{eqnarray}
where the sum is over all the equivalence classes $[(W,Q)]$  of the $r$-relevant pairs where $Q$ has $\alpha(Q)\leq N$ parts.
\end{corollary}


\bigskip
\section{Counting planar graphs (Proof of Theorem \ref{thm.main1})}
\label{sec.planargraph}
In this section we prove Theorem \ref{thm.main1}. To this end we first derive a relation (Theorem 
\ref{thm.main9}) between  the enumeration of planar graphs and the Gaussian integration of the generating function for edge-disjoint closed trails. We then (in Theorem \ref{thm.ice} and Corollary \ref{cor.labeledplanar}) derive an identity between this Gaussian integration and the Gaussian integration of the ice-type partition function for eulerian sets defined in (\ref{def.eta}).  Theorem \ref{thm.main1} follows from Corollary \ref{cor.labeledplanar}.

\bigskip
We consider the generating function for edge-disjoint closed trails definded by
\begin{eqnarray}
\zeta(M,x,y):= \sum_{r\ge 0}\omega_r(M) {x^r} \label{def.zeta}
\end{eqnarray}
where $\omega_r(M)\; = \; \sum_{q\in {\mathcal W}(r)}y^{|q|/2}\prod_{e\in q}M_e$ (as defined in (\ref{def.combfunc})) and ${\mathcal W}(r)$ is the set of all subsets of edges of $D(M)$ that can be decomposed into 
$r$ edge-disjoint closed trails.

\begin{theorem}
\label{thm.main9}
Let $s(N)\le \frac{N}{2}$ be any function in $N$ which approaches to $\infty$ slower than $N$ when $N\to \infty$.
The Gaussian matrix integral of the function $\zeta$ defined in (\ref{def.zeta}) satisfies the following.
\begin{eqnarray*}
&& \exp\left(\sum_{r\ge 0 } \sum_{e\ge 0 }{z^{r}y^{e-r}}\left(E_1(N,r,e)+\sum_{[G_{r,e}]}\frac{1}{|{\rm Aut}(G_{r,e})|}\right)\right) 
\ ' \leq_{s(N)}\, ' \ \\
&&e^{-N^2}\langle \zeta(M, Nzy^{-1},y) \rangle \\
&&' \leq_{\frac{N}{2}}\, ' \
\exp\left(\sum_{r\ge 0 } \sum_{e\ge 0 }{z^{r}y^{e-r}}\left(E_2(N,r,e)+\sum_{[G_{r,e}]}\frac{1}{|{\rm Aut}(G_{r,e})|}\right)\right),
\end{eqnarray*}
where the sum is over all isomorphism classes $[G_{r,e}]$ of connected simple graphs  which have planar embeddings with $r$ faces, $e$ edges and $\le N$ vertices, and $E_1(N,r,e), E_2(N,r,e)$ are of type $O(N^{-1})$. \end{theorem}

The following theorem formulates  the enumeration of planar graphs as the Gaussian integration of an ice-type partition function via eulerian sets.
\begin{theorem}
\label{thm.ice}
Let $\zeta$ be the function defined in (\ref{def.zeta}) and consider a function $\eta$ which maps $M\in \mathcal H_N$ to 
\begin{eqnarray}
\eta(M,x,y):= \sum_{A\subset N\times N \text{ eulerian}}U(A,x)y^{|A|/2}\prod_{e\in A}M_e,\label{def.icefunc}
\end{eqnarray}
where $U(A,x)= \prod_{C: \text{ component of } A}(\sum_{i=1}^{|C|/2}{x^i})$. Then
$$
\langle \zeta \rangle \; \; = \;  \; \langle \eta \rangle .
$$
\end{theorem}
\begin{proof}
This is true from the Euler theorem without $\langle \cdot \rangle $ when we restrict the eulerian subsets of $D$ to the subsets consisting of pairs of
oppositely directed edges with the same endvertices. In the integration the remaining terms of both $\zeta,\eta$
disappear.
\end{proof}

Note that the function $\eta(M,x,y)$ defined in (\ref{def.icefunc})
is equal to
\begin{eqnarray*}
&& \sum_{A\subset N\times N \text{ eulerian}} (\prod_{e\in A}M_e) \prod_{C: \text{ component of } A} \frac{x(x^{|C|/2}-1)}{x-1}y^{|C|/2}.
\end{eqnarray*}
Thus we have
\begin{eqnarray*}
\eta(M,Nzy^{-1},y) 
&=& \sum_{A\subset N\times N \text{ eulerian}} (\prod_{e\in A}M_e) \prod_{C: \text{ component of } A} \frac{(Nz)^{|C|/2+1}-(Nz)y^{|C|/2}}{Nz-y},
\end{eqnarray*}
as defined in (\ref{def.eta}). As a consequence of this and Theorem \ref{thm.ice} we have the following.
\begin{corollary}\label{cor.labeledplanar}
Theorem \ref{thm.main9} is true when we replace $\zeta(M, Nzy^{-1},y)$ by $\eta(M, Nzy^{-1},y)$.
\end{corollary}

\bigskip
In the rest of this section we will prove Theorem \ref{thm.main9} and then Theorem \ref{thm.main1}.

\begin{proof} (of Theorem \ref{thm.main9}) From Corollary ~\ref{cor.main3} we get
\begin{eqnarray}
&&\langle \zeta(M, Nzy^{-1},y) \rangle \nonumber\\
&=&\sum_{r\ge 0}\langle \omega_r \rangle  {(Nzy^{-1})^{r}} \nonumber\\
&\stacksign{(\ref{eq.equivalent})}{=}&
\sum_{r\ge 0}\sum_{[(W,Q)]} N^{r-e(W)}\frac{N(N-1)\dots (N-\alpha(Q)+1)}{|{\rm Aut}(G_Q(W^*))|}{z^{r}y^{e(W)-r}} \nonumber\\
&=&
\sum_{r\ge 0}\sum_{[(W,Q)]}N^{r-e(W)}N^{\alpha(Q)}\prod_{1\le i \le \alpha(Q)-1}\left(1-\frac{i}{N}\right)
 \frac{z^{r}y^{e(W)-r}}{|{\rm Aut}(G_Q(W^*))|}, \label{eq.rell}
\end{eqnarray}
where the sum is over all equivalence classes $[(W,Q)]$ of the $r$-relevant pairs where $Q$ has $\alpha(Q)\le N$ parts.

To see how (\ref{eq.rell}) gets non-zero contributions from equivalence classes of the $r$-relevant pairs, observe the following. Each equivalence class $[(W,Q)]$ can be represented by a fat graph $W$ with some  number $f(W)$ of faces satisfying  $\alpha(Q) \le f(W)\le 2e(W)$. 

For further computation we take any representation $(W,Q)_R$ of each $[(W,Q)]$ and get
\begin{eqnarray}
&&\langle \zeta(M, Nzy^{-1},y) \rangle \nonumber\\
&&=\sum_{r\ge 0}\sum_{(W,Q)_R}N^{r-e(W)+f(W)}
N^{\alpha(Q)-f(W)}\prod_{1\le i \le \alpha(Q)-1}\left(1-\frac{i}{N}\right)\frac{z^{r}y^{e(W)-r}}{|{\rm Aut}(G_Q(W^*))|}.\label{eq.rellrep}
\end{eqnarray}

We will bound  (\ref{eq.rellrep}) from above and from below.

Let $S_1(N,y,z)$ and $S_2(N,y,z)$ be defined as follows:
\begin{eqnarray}
&&S_1(N,y,z)
=\sum_{r\ge 0}\sum_{(W,Q)_R}N^{r-e(W)+f(W)}N^{\alpha(Q)-f(W)}
\left(1-\frac{2s(N)}{N}\right)^{\alpha(Q)}\frac{z^{r}y^{e(W)-r}}{|{\rm Aut}(G_Q(W^*))|}\nonumber\\ 
&&S_2(N,y,z)
=\sum_{r\ge 0}\sum_{(W,Q)_R}N^{r-e(W)+f(W)}N^{\alpha(Q)-f(W)}\frac{z^{r}y^{e(W)-r}}{|{\rm Aut}(G_Q(W^*))|}
\end{eqnarray}
where $s(N)$ is any function in $N$ which approaches to $\infty$ slower than $N$ when $N\to \infty$.

It is straightforward (we recall that that $\alpha(Q)\leq 2e(W)$) that
\begin{eqnarray*}
S_1(N,y,z)\ '\leq_{s(N)}\, '\  \langle \zeta(M, Nzy^{-1},y) \rangle \ '\leq\, '\  S_2(N,y,z).
\end{eqnarray*}

First we study $ S_2(N,y,z)$.
Recall that for each representation $(W,Q)_R$, $W$ is a fat graph with $r$ labelled vertices and $f(W)$  faces. Define the set $S_{(W,Q)_R}$ of siblings $(W',Q)_R$ of each $(W,Q)_R$, where $W'$ is obtained from $W$ by relabelling of its vertex labels. Since $\sum_{W'\in S_{(W,Q)_R}}1=r!/|{\rm Aut}(W)|$ (and notice that $v(W)=v(W')=r, e(W)=e(W'), f(W)=f(W')$, $|{\rm Aut}(W)|=|{\rm Aut}(W')|$ and $|{\rm Aut}(G_{Q}(W^*))|=|{\rm Aut}(G_{Q}(W'^*))|$), we have that
\begin{eqnarray}
&&S_2(N,y,z)=\sum_{r\ge 0}\sum_{(W,Q)_R} \left(\frac{\sum_{W'\in S_{(W,Q)_R}} 1}{r!/ |{\rm Aut}(W)|} \right) N^{r-e(W)+ f(W)}
N^{\alpha(Q)- f(W)}\frac{z^{r}y^{e(W)-r}}{|{\rm Aut}(G_{Q}(W^*))|}\nonumber\\
&&
=\sum_{r\ge 0}\sum_{(W,Q)_R}\; \; \sum_{W'\in S_{(W,Q)_R}}N^{r-e(W')+ f(W')}
N^{\alpha(Q)- f(W')}\frac{|{\rm Aut}(W')|}{|{\rm Aut}(G_{Q}(W'^*))|} \frac{z^{r}y^{e(W')-r}}{r!}.\label{eq.towardconnectedlsibling}
\end{eqnarray}

Let $\{G_1,\cdots,G_c\}$ be the set of the connected components of $G_Q(W^*)$. This induces a partition of the components of $W$ into $W_1,\ldots,W_c$,  and a partition of the classes of $Q$ into $Q_1,\ldots,Q_c$. Each $Q_i$ is a partition of the faces of $W_i$.

Notice in (\ref{eq.towardconnectedlsibling}) that each of the factors,
$N^{r-e(W')+ f(W')}N^{\alpha(Q)- f(W')}=N^{v(W')-e(W')+ \alpha(Q)}$, 
$|{\rm Aut}(W')|$, $|{\rm Aut}(G_{Q}(W'^*))|$, $z^ry^{e(W')-r}=z^{v(W')}y^{e(W')-v(W')}$  
is  multiplicative under disjoint union of $(W_i,Q_i)$.
We have 
${r!}/{|{\rm Aut}(W)|}$ 
ways to relabelling the vertices of $W$, which is equal to 
$$\binom{r}{r_1, r_2,\cdots,r_c}\times \frac{r_1!\times  r_2!\times\cdots\times  r_c!}{|{\rm Aut}(W_1)| \times  |{\rm Aut}(W_2)| \times \cdots \times |{\rm Aut}(W_c)|},$$ 
that is, the number of ways of choosing  vertex sets for each $W_i$ and labelling $W_i$. Thus we can use the relation between the generating function for not-necessarily connected representations and the exponent of the generating function for connected ones to obtain
\begin{eqnarray}
&&S_2(N,y,z)=_{\frac{N}{2}}\exp(T_2(N,y,z)) \label{eq.S2T2}
 \end{eqnarray}
with
\begin{eqnarray*}
&&T_2(N,y,z)=\sum_{r\ge 0} \sum_{(W,Q)_R^c}\sum_{W'\in S_{(W,Q)_R^c}}N^{r-e(W')+ f(W')}N^{\alpha- f(W')}\frac{|{\rm Aut}(W')|}{|{\rm Aut}(G_Q(W'^*))|} \frac{z^{r}y^{e(W')-r}}{r!} 
 \end{eqnarray*}
where the second sum is over all representations $(W,Q)_R^c$ where $G_Q(W^*)$ is {\em connected}. The bound $\frac{N}{2}$ in (\ref{eq.S2T2}) is imposed, since in $S_2$ the total degree in $y$ and $z$ (which equals $e(W')=e(W)$) $\le\frac{N}{2} $ implies that $\textcolor{black}{f}(W)\le N$. Observe that for any nimple graph, $\textcolor{black}{f}(W)\le 2e(W)$.

\smallskip
Since $\sum_{W'\in S_{(W,Q)_R^c}} 1 = {r!}/{|{\rm Aut}(W')|} $,  we have
\begin{eqnarray} 
&&T_2(N,y,z)\nonumber\\ 
&&=\sum_{r\ge 0} \quad \sum_{(W,Q)_R^c }N^{r-e(W)+ f(W)}N^{\alpha(Q)- f(W)}\frac{z^{r}y^{e(W)-r}}{|{\rm Aut}(G_Q(W^*))|}\sum_{W'\in S_{(W,Q)_R^c}}\frac{|{\rm Aut}(W')|}{r!}\nonumber\\
&&=\sum_{r\ge 0}\quad  \sum_{(W,Q)_R^c}N^{r-e(W)+ f(W)}
N^{\alpha(Q)- f(W)}
 \frac{z^{r}y^{e(W)-r}}{|{\rm Aut}(G_Q(W^*))|}\nonumber\\
&&=\sum_{r\ge 0} \sum_{e\ge 0}{z^{r}y^{e-r}}\quad \sum_{(W,Q)_R^c }N^{2-2g(W)}
N^{\alpha(Q)- f(W)}
 \frac{1}{|{\rm Aut}(G_Q(W^*))|},\label{eq.T2}
 \end{eqnarray}
where  $W$ has $r$ vertices, $e$ edges, and $\le N$ faces.

We recall that for each $r\ge 0$,
$$
2-2g(W)+ \alpha(Q) -f(W)= v(W)- e(W)+ \alpha(Q)= r-e(G_Q(W^*))+v(G_Q(W^*))
$$
is independent of the choice of the representation $(W,Q)_R$.

If $(W,Q)_R^c$ has $g(W)>0$, or $g(W)=0$ and $Q$ is non-trivial,
then 
$$N^{-2g(W)}N^{\alpha(Q)- f(W)}=O(N^{-1}).$$
The contribution of such representations to $T_2(N,y,z)$ is  
\begin{eqnarray} 
&&N^2 \sum_{r\ge 0} \sum_{e\ge 0}{z^{r}y^{e-r}}\quad E_2(N,r,e),  \label{eq.genusrep1} 
\end{eqnarray}
where 
\begin{eqnarray} 
E_2(N,r,e)
&&=\sum\nolimits_{[(W,Q)^c] }^{\text{err}}O(N^{-1}) \frac{1}{|{\rm Aut}(G_Q(W^*))|},
\label{eq.E2}
\end{eqnarray}
where the sum is over all the equivalence classes $[(W,Q)^c]$  of the $r$-relevant pairs where $G_Q(W^*)$ is {connected},  $W$ has $e$ edges, and $g(W)>0$ or $Q$ non-trivial.

\smallskip
On the other hand, for the remaining  representations $(W,Q)_R^c$ (that is, $g(W)=0$ and  $Q$ is trivial), 
$$
N^{-2g(W)}N^{\alpha(Q)- f(W)}=1.
$$
Hence the contribution of such representations to $T_2(N,y,z)$ is
\begin{eqnarray}
&&N^2\sum_{r\ge 0} \sum_{e\ge 0}{z^{r}y^{e-r}}\sum_{[(  W,\emptyset)^c]\, \textit{planar}}\frac{1}{|{\rm Aut}(W^*)|}.\label{eq.planarcontr}
 \end{eqnarray}

We can also sum over the duals in (\ref{eq.planarcontr}). The dual of a {\em connected} fat graph is {\em connected} and the dual of a {\em relevant} fat graph is {\em nimple}. Moreover, there is exactly one nimple graph on at most one vertex  (with no edge) and at most one simple connected graph with one vertex and no edge. Hence (\ref{eq.planarcontr})  is equal to
 \begin{eqnarray}
N^2\sum_{r\ge 0}\sum_{e\ge 0}{z^{r}y^{e-r}}\sum_{[G_{r,e}]} \frac{1}{|{\rm Aut}(G_{r,e})|},\label{eq.planarcontr1}
 \end{eqnarray}
where the sum is over all isomorphism classes $[G_{r,e}]$ of {\em connected simple} graphs which have planar embedding (we recall that $Q$ is trivial) with $r$ faces, $e$ edges and $\le N$ vertices.

Putting (\ref{eq.S2T2}), (\ref{eq.genusrep1}), and (\ref{eq.planarcontr1}) together, we obtain that
 \begin{eqnarray}
S_2(N,y,z)=_{\frac{N}{2}}\exp\left(N^2 \sum_{r\ge 0}\sum_{e\ge 0}{z^{r}y^{e-r}}\left(E_2(N,r,e)+\sum_{[G_{r,e}]} \frac{1}{|{\rm Aut}(G_{r,e})|}\right)\right).\label{eq.upper}
 \end{eqnarray}

The analysis for $S_1(N,y,z)$ is analogous as to $S_2(N,y,z)$; we only note that the factor $\left(1-\frac{2s(N)}{N}\right)^{\alpha(Q)}$ is multiplicative. We obtain that
\begin{eqnarray*}
S_1(N,y,z)=_{\frac{N}{2}}\exp\left(N^2 \sum_{r\ge 0}\sum_{e\ge 0}{z^{r}y^{e-r}}\left(E_1(N,r,e)+\sum_{[G_{r,e}]} \frac{1}{|{\rm Aut}(G_{r,e})|}\right)\right),\label{eq.lower}
 \end{eqnarray*}
where
\begin{eqnarray} 
E_1(N,r,e)
&&=\sum\nolimits_{[(W,Q)^c] }^{{\text{err}}}O(N^{-1}) \left(1-\frac{2s(N)}{N}\right)^{\alpha(Q)}\frac{1}{|{\rm Aut}(G_Q(W^*))|}\nonumber\\
&&+\sum_{[G_{r,e}]} \left(-1+\left(1-\frac{2s(N)}{N}\right)^{\alpha(Q)}\right)\frac{1}{|{\rm Aut}(G_{r,e})|},\label{eq.E1}
\end{eqnarray}
where the sum in the first term is as in $E_2(N,r,e)$, and the sum in the second term is as in (\ref{eq.planarcontr1}).
\end{proof}

\smallskip

\begin{proof} (of Theorem \ref{thm.main1}) From Corollary \ref{cor.labeledplanar} we have that
\begin{eqnarray}
&& \exp\left(\sum_{r\ge 0 } \sum_{e\ge 0 }{z^{r}y^{e-r}}\left(E_1(N,r,e)+\sum_{[G_{r,e}]}\frac{1}{|{\rm Aut}(G_{r,e})|}\right)\right) 
\ ' \leq_{s(N)}\, ' \ \nonumber\\
&&e^{-N^2}\langle \eta(M, Nzy^{-1},y) \rangle \nonumber\\
&&' \leq_{\frac{N}{2}}\, ' \
\exp\left(\sum_{r\ge 0 } \sum_{e\ge 0 }{z^{r}y^{e-r}}\left(E_2(N,r,e)+\sum_{[G_{r,e}]}\frac{1}{|{\rm Aut}(G_{r,e})|}\right)\right).\label{eq.Gedge}
\end{eqnarray}

We first study the leading term $\sum_{[G_{r,e}]}{1}/{|{\rm Aut}(G_{r,e})|}$ in (\ref{eq.Gedge}). We note that to each isomorphism class $[G_{r,e}]$ of connected simple graphs which have planar embeddings with $r$ faces, $e$ edges, and $v(G_{r,e})\le N$ vertices, there are $v(G_{r,e})!\ /\ |{\rm Aut}(G_{r,e})|$ many labelled corresponding ones. From this, together with the Euler relation, $v(G_{r,e})-e+r=2$, we have that 
\begin{eqnarray}
&&\sum_{r\ge 0 } \sum_{e\ge 0 }{z^{r}y^{e-r}} \sum_{[G_{r,e}]}\frac{1}{|{\rm Aut}(G_{r,e})|}\nonumber\\
&&=\sum_{r\ge 0 }\sum_{e\ge 0 }\quad \sum_{G_{r,e}\, \,  \text{labelled}} \frac{z^{r}y^{v(G_{r,e})-2}}{|{\rm Aut}(G_{r,e})|} \times \frac{|{\rm Aut}(G_{r,e})|}{v(G_{r,e})!}\nonumber\\
&&=\sum_{r\ge 0 }\sum_{e\ge 0 }\quad \sum_{G_{r,e}\, \,  \text{labelled}} \frac{z^{r}y^{v(G_{r,e})-2}} {v(G_{r,e})!}. \label{eq.labelled}
\end{eqnarray}
Being sorted by a possible number $n$ of vertices,  (\ref{eq.labelled}) is equal to
\begin{eqnarray*}
\sum_{n= 1}^{N} \sum_{r\ge 0 }p(n,r)\frac{z^{r}y^{n-2}}{n!},
\end{eqnarray*}
where $p(n,r)$ denotes the number of labelled connected simple graphs on $n$ vertices that have planar embeddings with $r$ faces.

Concerning the subleading term $E_2(N,r,e)$ in (\ref{eq.Gedge}), 
due to the Euler relation, $e-r=f(W)+2g(W)-2$, we have
\begin{eqnarray}
&&\sum_{r\ge 0 } \sum_{e\ge 0 }{z^{r}y^{e-r}}E_2(N,r,e)\nonumber\\
&&\stacksign{(\ref{eq.E2})}{=}\; \sum_{r\ge 0 } \sum_{e\ge 0 }{z^{r}y^{e-r}}\sum\nolimits_{[(W,Q)^c]}^{\text{err}}\quad \frac{O(N^{-1})}{|{\rm Aut}(G_Q(W^*))|}\nonumber\\
&&\; = \; \sum_{r\ge 0 }{z^{r}y^{f(W)+2g(W)-2}}\sum\nolimits_{[(W,Q)^c] }^{\text{err}}\frac{O(N^{-1}) }{|{\rm Aut}(G_Q(W^*))|}.\label{eq.GE2}
\end{eqnarray}
If (\ref{eq.GE2}) is sorted by a possible number $n$ where $n=f(W)+2g(W)$, it is equal to
\begin{eqnarray*}
&&\sum_{n\ge 1} \sum_{r\ge 0 }e_2(N,n,r)\frac{z^{r}y^{n-2}}{n!},
\end{eqnarray*}
with $$e_2(N,n,r):=\sum\nolimits_{[(W,Q)^c]}^{\text{err}} O(N^{-1}) \frac{n!}{|{\rm Aut}(G_Q(W^*))|},$$ 
where the sum is over all equivalent classes $[(W,Q)^c]$ of the $r$-relevant pairs where $G_Q(W^*)$ is connected, $W$ satisfies $n=f(W)+2g(W)$,  and $g(W)>0$ or $Q$ non-trivial. Note that the sum is finite, since $e(W)=n+r-2$.

The analysis of $E_1(N,r,e)$ and $e_1(N,n,r)$ is analogous.
\end{proof}


\bigskip
\section{Counting directed cycle double covers (Proof of Theorem \ref{thm.main2})}
\label{sec.dcdc}
In this section we prove Theorem \ref{thm.main2}. To this end we first formulate the enumeration of graphs with a {\em specified} directed cycle double cover (DCDC) as the Gaussian matrix integral of the generating function $\xi$ for even sets with its specified  cycle decomposition (defined in (\ref{fun.dcdc})). We further show that $\xi$ is indeed an Ihara-Selberg type function,  via a lemma on coin arrangements (Lemma \ref{lem:coin}).

First we recall from the introduction that 
a subset $A$ of directed edges of $D$ is {\it even} if $A$  can be written as a union of edge-disjoint directed cycles  of length bigger than two. 
We further denote by ${\mathcal A}$ the set of all pairs $(q,K)$ where $q$ is an even set of edges of $D$
and $K$ is a decomposition of $q$ into directed cycles of length at least three.
The following observation is straightforward.

\begin{observation}
\label{l.dc}
Let $(q,K)\in {\mathcal A}$ and let $P$ be a proper pairing of $q$. Then the cycles of $K$ form a DCDC of a nimple graph with vertices $\{1,\dots, N\}$  and the edges given by the pairing $P$.
\end{observation}

\begin{definition}
\label{def.we}
Let $G$ be a finite nimple graph with at most $N$ vertices and let $C$ be a DCDC of $G$. 
Then let $c(G,C)$ be the set of all triples $(q,P,C)$ so that there is a colouring $d$ of the vertices 
of $G$ by colours $\{1,\dots, N\}$, where each vertex gets a different colour, $$q= \{(d(x),d(y)),(d(y),d(x));\{x,y\}\in E(G)\},$$
and $P$ consists of all the pairs $([(d(x),d(y)),(d(y),d(x))];\{x,y\}\in E(G))$.
\end{definition}

We remark that each such pair $(q,C)$ in Definition \ref{def.we} belongs to ${\mathcal A}$
and that $$|c(G,C)|\; = \; \frac{N(N-1)\cdots (N-|V(G)|+1)}{|{\rm Aut}(G,C)|}.$$ 

\begin{proposition}
\label{prop.la}
For the function $\xi$ defined in (\ref{fun.dcdc}) which maps $M\in \mathcal H_N$ to
$$\xi(M):= \; \sum_{(q,C)\in {\mathcal A}}\prod_{e\in q}M_e,$$ 
a term $((q,C),P)$ with $P$ a proper pairing of $q$ contributes to $\langle \xi \rangle $ if and only if
 there is a nimple graph $G$ with at most $N$ vertices where $C$ is its DCDC and such that $(q,P,C) \in c(G,C)$.
\end{proposition}
\begin{proof}
If $(q,P,C)\in c(G,C)$, then $C$ provides a partition of $q$ into its directed cycles 
and hence $((q,C),P)$ contributes to $\langle \xi \rangle $. 
On the other hand if $((q,C),P)$ contributes to $\langle \xi \rangle $, then letting $G$ be the graph with the vertices from $\{1,\dots, N\}$ and the edges given by $P$ we get that $G$ is nimple. Moreover $q$ comes with a decomposition $C$ into the 
directed cycles and thus $(q,P,C)\in c(G,C)$.
  \end{proof}

 \begin{proposition}
\label{prop.lla}
If $c(G,C)\cap c(G',C')\neq \emptyset$, then $(G,C)$ is isomorphic to $(G',C')$. 
Moreover, if $(G,C)$ is isomorphic to $(G',C')$, then $c(G,C)= c(G',C')$.
\end{proposition}
\begin{proof}
If $(q,P,C)\in c(G,C)\cap c(G',C')$, then the construction of $q$ induces a function between the sets of vertices of 
$G$ and $G'$, and $P$ gives the edges of both $G, G'$. Hence there is an isomorphism from $G$ to $G'$ which preserves
the fixed DCDCs.
The second part is true since the definition of $c(G,C)$ does not depend on 'names' of the vertices. 
  \end{proof}

As a consequence we have the following.

\begin{theorem}
\label{thm.main4}
Let $\xi$ be the function defined in (\ref{fun.dcdc}). Then
$$
\langle \xi \rangle \; = \; \sum_{[(G,C)]}\frac{N(N-1)\dots (N-|V(G)|+1)}{|{\rm Aut}(G,C)|N^{e(G)}},
$$
where the sum is over all isomorphism classes of pairs $(G,C)$ where $G$ is a nimple graph with at most $N$ vertices and $C$ {\em a specified  DCDC} of $G$.
\end{theorem}

\bigskip
\subsection{Calculations}
\label{subsec.ihara}
The integral $\langle \xi \rangle $ counts all the directed cycle double covers of 
graphs on at most $N$ vertices and hence its calculation is an attractive task which need not be
hopeless. We show next a curious formula for $\xi(M)$ which identifies it with an
Ihara-Selberg-type function  (see Theorem \ref{thm.prr}).

\bigskip
{\bf Construction of digraph $D'$.} We first construct a directed graph $D'$ with the weights on
 the {\it transitions between the edges}.

\begin{figure}[h]
\begin{center}
\begin{picture}(0,0)%
\includegraphics{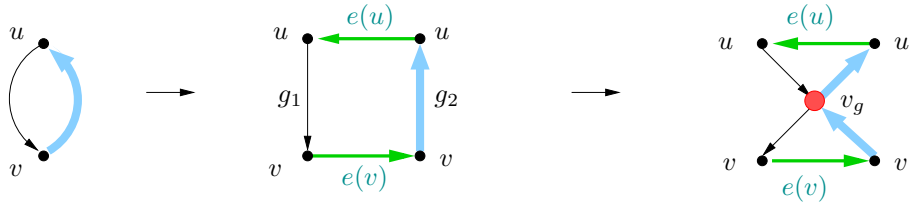}%
\end{picture}%
\setlength{\unitlength}{2693sp}%
\begingroup\makeatletter\ifx\SetFigFont\undefined%
\gdef\SetFigFont#1#2#3#4#5{%
  \reset@font\fontsize{#1}{#2pt}%
  \fontfamily{#3}\fontseries{#4}\fontshape{#5}%
  \selectfont}%
\fi\endgroup%
\begin{picture}(8524,1806)(2376,-3994)
\put(2386,-3751){\makebox(0,0)[lb]{\smash{{\SetFigFont{10}{12.0}{\familydefault}{\mddefault}{\updefault}{\color[rgb]{0,0,0}$v$}%
}}}}
\put(6301,-2491){\makebox(0,0)[lb]{\smash{{\SetFigFont{10}{12.0}{\familydefault}{\mddefault}{\updefault}{\color[rgb]{0,0,0}$u$}%
}}}}
\put(6301,-3076){\makebox(0,0)[lb]{\smash{{\SetFigFont{10}{12.0}{\familydefault}{\mddefault}{\updefault}{\color[rgb]{0,0,0}$g_2$}%
}}}}
\put(10036,-3121){\makebox(0,0)[lb]{\smash{{\SetFigFont{10}{12.0}{\familydefault}{\mddefault}{\updefault}{\color[rgb]{0,0,0}$v_g$}%
}}}}
\put(6346,-3661){\makebox(0,0)[lb]{\smash{{\SetFigFont{10}{12.0}{\familydefault}{\mddefault}{\updefault}{\color[rgb]{0,0,0}$v$}%
}}}}
\put(8956,-3706){\makebox(0,0)[lb]{\smash{{\SetFigFont{10}{12.0}{\familydefault}{\mddefault}{\updefault}{\color[rgb]{0,0,0}$v$}%
}}}}
\put(8911,-2581){\makebox(0,0)[lb]{\smash{{\SetFigFont{10}{12.0}{\familydefault}{\mddefault}{\updefault}{\color[rgb]{0,0,0}$u$}%
}}}}
\put(2386,-2491){\makebox(0,0)[lb]{\smash{{\SetFigFont{10}{12.0}{\familydefault}{\mddefault}{\updefault}{\color[rgb]{0,0,0}$u$}%
}}}}
\put(4816,-2491){\makebox(0,0)[lb]{\smash{{\SetFigFont{10}{12.0}{\familydefault}{\mddefault}{\updefault}{\color[rgb]{0,0,0}$u$}%
}}}}
\put(4771,-3751){\makebox(0,0)[lb]{\smash{{\SetFigFont{10}{12.0}{\familydefault}{\mddefault}{\updefault}{\color[rgb]{0,0,0}$v$}%
}}}}
\put(10531,-2581){\makebox(0,0)[lb]{\smash{{\SetFigFont{10}{12.0}{\familydefault}{\mddefault}{\updefault}{\color[rgb]{0,0,0}$u$}%
}}}}
\put(10531,-3706){\makebox(0,0)[lb]{\smash{{\SetFigFont{10}{12.0}{\familydefault}{\mddefault}{\updefault}{\color[rgb]{0,0,0}$v$}%
}}}}
\put(5491,-2356){\makebox(0,0)[lb]{\smash{{\SetFigFont{10}{12.0}{\familydefault}{\mddefault}{\updefault}{\color[rgb]{0,.56,.56}$e(u)$}%
}}}}
\put(5446,-3841){\makebox(0,0)[lb]{\smash{{\SetFigFont{10}{12.0}{\familydefault}{\mddefault}{\updefault}{\color[rgb]{0,.56,.56}$e(v)$}%
}}}}
\put(9496,-3931){\makebox(0,0)[lb]{\smash{{\SetFigFont{10}{12.0}{\familydefault}{\mddefault}{\updefault}{\color[rgb]{0,.56,.56}$e(v)$}%
}}}}
\put(9541,-2356){\makebox(0,0)[lb]{\smash{{\SetFigFont{10}{12.0}{\familydefault}{\mddefault}{\updefault}{\color[rgb]{0,.56,.56}$e(u)$}%
}}}}
\put(4861,-3076){\makebox(0,0)[lb]{\smash{{\SetFigFont{10}{12.0}{\familydefault}{\mddefault}{\updefault}{\color[rgb]{0,0,0}$g_1$}%
}}}}
\end{picture}%
\end{center}
\caption{Construction of $e(v)$ and $v_g$}\label{fig:joint}
\end{figure}

First we split each vertex of $D$, i.e., we replicate each vertex $v$, and then we connect two $v$'s by a new edge $e(v)$
and we let all the edges of $D$ entering $v$ enter the initial vertex of $e(v)$, and
all the edges of $D$ leaving $v$ leave the terminal vertex of $e(v)$.
If an edge $g_1$ enters $v$ in $D$ then we define the weight of the transition $w(g_1,e(v)) = M_{g_1}$.
We let all the remaining transition be equal to one (see Figure~\ref{fig:joint}, the first two parts).

Finally, for each pair $g_1,g_2$ of oppositely directed edges of $D$, say $g_1=(uv),g_2=(vu)$  
we introduce a new vertex $v_g$ (with $g=\{u,v\}$) and we let both $g_1,g_2$ pass through it; equivalently, we subdivide both $g_1,g_2$ by one vertex and identify this pair of vertices into a unique vertex called $v_g$ 
(and thus we have new edges $(uv_g),(v_gv)$ from $g_1=(uv)$, and new edges $(vv_g),(v_gu)$ from $g_2=(vu)$) (see Figure~\ref{fig:joint}, the last two parts). 

\begin{figure}[h]
\begin{center}
\begin{picture}(0,0)%
\includegraphics{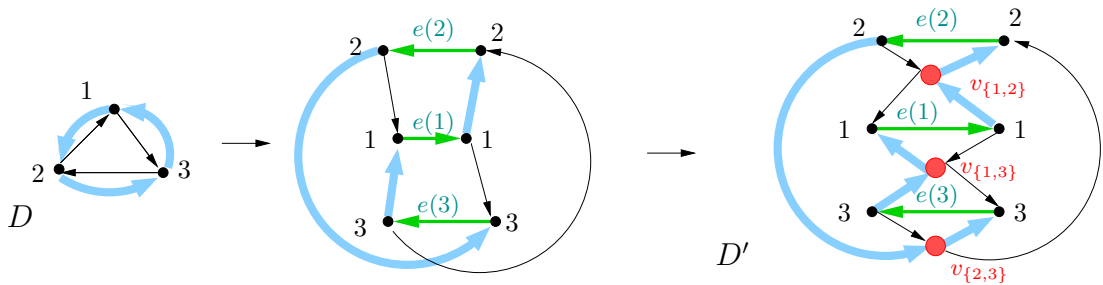}%
\end{picture}%
\setlength{\unitlength}{2693sp}%
\begingroup\makeatletter\ifx\SetFigFont\undefined%
\gdef\SetFigFont#1#2#3#4#5{%
  \reset@font\fontsize{#1}{#2pt}%
  \fontfamily{#3}\fontseries{#4}\fontshape{#5}%
  \selectfont}%
\fi\endgroup%
\begin{picture}(10075,2497)(1666,-3899)
\put(1666,-3391){\makebox(0,0)[lb]{\smash{{\SetFigFont{12}{14.4}{\familydefault}{\mddefault}{\updefault}{\color[rgb]{0,0,0}$D$}%
}}}}
\put(3241,-2941){\makebox(0,0)[lb]{\smash{{\SetFigFont{10}{12.0}{\familydefault}{\mddefault}{\updefault}{\color[rgb]{0,0,0}$3$}%
}}}}
\put(5401,-1636){\makebox(0,0)[lb]{\smash{{\SetFigFont{9}{10.8}{\familydefault}{\mddefault}{\updefault}{\color[rgb]{0,.56,.56}$e(2)$}%
}}}}
\put(5401,-2491){\makebox(0,0)[lb]{\smash{{\SetFigFont{9}{10.8}{\familydefault}{\mddefault}{\updefault}{\color[rgb]{0,.56,.56}$e(1)$}%
}}}}
\put(5446,-3256){\makebox(0,0)[lb]{\smash{{\SetFigFont{9}{10.8}{\familydefault}{\mddefault}{\updefault}{\color[rgb]{0,.56,.56}$e(3)$}%
}}}}
\put(10036,-1546){\makebox(0,0)[lb]{\smash{{\SetFigFont{9}{10.8}{\familydefault}{\mddefault}{\updefault}{\color[rgb]{0,.56,.56}$e(2)$}%
}}}}
\put(9856,-2401){\makebox(0,0)[lb]{\smash{{\SetFigFont{9}{10.8}{\familydefault}{\mddefault}{\updefault}{\color[rgb]{0,.56,.56}$e(1)$}%
}}}}
\put(10036,-3166){\makebox(0,0)[lb]{\smash{{\SetFigFont{9}{10.8}{\familydefault}{\mddefault}{\updefault}{\color[rgb]{0,.56,.56}$e(3)$}%
}}}}
\put(10351,-3841){\makebox(0,0)[lb]{\smash{{\SetFigFont{9}{10.8}{\familydefault}{\mddefault}{\updefault}{\color[rgb]{1,0,0}$v_{\{2,3\}}$}%
}}}}
\put(10441,-2896){\makebox(0,0)[lb]{\smash{{\SetFigFont{9}{10.8}{\familydefault}{\mddefault}{\updefault}{\color[rgb]{1,0,0}$v_{\{1,3\}}$}%
}}}}
\put(10531,-2131){\makebox(0,0)[lb]{\smash{{\SetFigFont{9}{10.8}{\familydefault}{\mddefault}{\updefault}{\color[rgb]{1,0,0}$v_{\{1,2\}}$}%
}}}}
\put(10936,-2581){\makebox(0,0)[lb]{\smash{{\SetFigFont{10}{12.0}{\familydefault}{\mddefault}{\updefault}{\color[rgb]{0,0,0}$1$}%
}}}}
\put(10936,-3346){\makebox(0,0)[lb]{\smash{{\SetFigFont{10}{12.0}{\familydefault}{\mddefault}{\updefault}{\color[rgb]{0,0,0}$3$}%
}}}}
\put(6121,-1681){\makebox(0,0)[lb]{\smash{{\SetFigFont{10}{12.0}{\familydefault}{\mddefault}{\updefault}{\color[rgb]{0,0,0}$2$}%
}}}}
\put(10891,-1546){\makebox(0,0)[lb]{\smash{{\SetFigFont{10}{12.0}{\familydefault}{\mddefault}{\updefault}{\color[rgb]{0,0,0}$2$}%
}}}}
\put(2341,-2221){\makebox(0,0)[lb]{\smash{{\SetFigFont{10}{12.0}{\familydefault}{\mddefault}{\updefault}{\color[rgb]{0,0,0}$1$}%
}}}}
\put(9451,-1636){\makebox(0,0)[lb]{\smash{{\SetFigFont{10}{12.0}{\familydefault}{\mddefault}{\updefault}{\color[rgb]{0,0,0}$2$}%
}}}}
\put(9316,-3346){\makebox(0,0)[lb]{\smash{{\SetFigFont{10}{12.0}{\familydefault}{\mddefault}{\updefault}{\color[rgb]{0,0,0}$3$}%
}}}}
\put(9316,-2581){\makebox(0,0)[lb]{\smash{{\SetFigFont{10}{12.0}{\familydefault}{\mddefault}{\updefault}{\color[rgb]{0,0,0}$1$}%
}}}}
\put(4816,-1816){\makebox(0,0)[lb]{\smash{{\SetFigFont{10}{12.0}{\familydefault}{\mddefault}{\updefault}{\color[rgb]{0,0,0}$2$}%
}}}}
\put(4951,-2671){\makebox(0,0)[lb]{\smash{{\SetFigFont{10}{12.0}{\familydefault}{\mddefault}{\updefault}{\color[rgb]{0,0,0}$1$}%
}}}}
\put(4861,-3481){\makebox(0,0)[lb]{\smash{{\SetFigFont{10}{12.0}{\familydefault}{\mddefault}{\updefault}{\color[rgb]{0,0,0}$3$}%
}}}}
\put(6256,-3436){\makebox(0,0)[lb]{\smash{{\SetFigFont{10}{12.0}{\familydefault}{\mddefault}{\updefault}{\color[rgb]{0,0,0}$3$}%
}}}}
\put(6031,-2671){\makebox(0,0)[lb]{\smash{{\SetFigFont{10}{12.0}{\familydefault}{\mddefault}{\updefault}{\color[rgb]{0,0,0}$1$}%
}}}}
\put(1891,-2986){\makebox(0,0)[lb]{\smash{{\SetFigFont{10}{12.0}{\familydefault}{\mddefault}{\updefault}{\color[rgb]{0,0,0}$2$}%
}}}}
\put(8191,-3751){\makebox(0,0)[lb]{\smash{{\SetFigFont{12}{14.4}{\familydefault}{\mddefault}{\updefault}{\color[rgb]{0,0,0}$D'$}%
}}}}
\end{picture}%
\end{center}
\caption{An example of the construction of $D'$ from $D$}\label{fig:const}
\end{figure}

We let the weights of the transitions at the vertex $v_g$ {\it between} $g_1$ and $g_2$ (i.e., between $(uv_g)$ and $(v_gu)$ and between $(vv_g)$ and $(v_gv)$) be equal to zero, 
the transitions {\it along} $g_1$ and $g_2$ (i.e., between $(uv_g)$ and $(v_gv)$ and between $(vv_g)$ and $(v_gu)$) be equal to one, 
and the transitions between $(v_gv)$ and $e(v)$  be equal to  $M_{g_1}$ and 
between $(v_gu)$ and $e(u)$  be equal to  $M_{g_2}$. 
See an example in Figure~\ref{fig:const}.

In what follows, the directed closed walk is considered {\it not pointed}. We let the weight
of the directed closed walk be the product of the weights of its transitions.

\begin{observation}
\label{o.d}
There is a weight preserving bijection between the set of the directed cycles of $D$
of length at least three and a non-zero weight, and the set of the aperiodic directed closed walks of $D'$ of a non-zero 
weight which go through
each directed edge and through each vertex $v_g$ at most once. 
\end{observation}
\begin{proof}
This follows directly from the construction of $D'$.
  \end{proof}

\begin{definition}
\label{def.rot}
We define the {\it rotation number} $\rho(w)$ for each closed walk $w$ of $D'$ with a non-zero weight
by induction as follows: first order the directed edges of $D'$, say as $a_1,\dots, a_m$,
so that the edges $e(v), v\in V(D)$ form the terminal segment. Then
\begin{enumerate}
\item[1.]
If $w$ is a directed cycle, then we let $\rho(w)= -1$.
\item[2.]
Let $w$ go at least twice through a directed edge. Let $a$ be the first such edge
in the fixed ordering.  Hence $w$ is a concatenation
of two shorter closed walks $w_1, w_2$, both containing $a$. If $a\neq e(v)$ for some $v$
then we let $\rho(w)= \rho(w_1)\rho(w_2)$. If $a= e(v)$, then we let $\rho(w)=0$.
\item[3.]
If none of 1.,2. applies, $w$ must go through the vertex $v_g$ (introduced in the definition
of $D'$) at least twice. Then we again let $\rho(w)= 0$.
\end{enumerate}
\end{definition}

\begin{theorem}
\label{thm.prr}
Let $\xi$ be the function defined in (\ref{fun.dcdc}). Then
$$
\xi(M)= \prod_p (1- \rho(p)w(p)),
$$
where the product is over all aperiodic directed closed walks $p$  in $D'$ and $w(p)$ denotes 
the weight of $p$.
\end{theorem}

To prove Theorem~\ref{thm.prr} we will need a curious lemma on coin arrangements stated below. It has been introduced
by Sherman \cite{S} in the study of 2-dimensional Ising problem.

\medskip\noindent

\begin{lemma}[A lemma on coin arrangements.]\label{lem:coin}
Suppose we have a fixed collection of $N$ objects of which $m_1$ are of
 one kind, $m_2$ are of second kind,  $\cdots$, and $m_n$ are of $n$-th kind. 
Let $b_k$ be the number of exhaustive unordered arrangements of these
symbols into $k$ disjoint, nonempty, circularly ordered sets such that
no two circular orders are the same and none are periodic. For example
let us have 10 coins of which 3 are pennies, 4 are nickles and 3 are quarters.
Then $\{(p,n),(n,p),(p,n,n,q,q,q)\}$ is not a correct arrangement since 
$(p,n)$ and $(n,p)$ represent the same circular order.
If $N>1$ then $\sum_{i=1}^N(-1)^{i+1}b_i=0$.
\end{lemma}

\smallskip
\begin{proof}[Proof of Lemma~\ref{lem:coin}]
The lemma follows immediately if we expand the LHS of the following \emph{Witt Identity} and collect
terms where the sums of the exponents of the $z_i$'s are the same. 

\bigskip
{\bf Witt Identity} (see \cite{H}): Let $z_1, ... ,z_k$ be commuting
variables. Then 
$$
\prod_{m_1, ... ,m_k\geq 0}(1-z_1^{m_1}...z_k^{m_k})^{M(m_1, ... ,m_k)}=1-z_1-z_2-...-z_k,
$$
where $M(m_1, .... ,m_k)$ is the number of different nonperiodic sequences of $z_i$'s
taken with respect to circular order.  
  \end{proof}

\smallskip
\begin{proof}[Proof of Theorem~\ref{thm.prr}]
We first show that the coefficients
corresponding to the products of variables where at least one $M_e, e\neq e(v)$, appears with the
exponent greater than one, are all equal to zero. 

Let us denote $W(p)= -\rho(p)w(p)$.
Let $A_1$ be the set of all aperiodic closed walks $p$ such that $a_1$
appears in $p$. Each $p\in A_1$ has a unique factorization into
words $(W_1,...,W_k)$ each of which starts with $a_1$ and has no other appearance of $a_1$.

Let $S$ be a monomial summand in the expansion of $\prod_{p\in A_1} (1+ W(p))$. 
Hence $S$ is a product of finitely many $W(p), p\in A_1$.

Each $p \in A_1$ has a unique factorization into words defined above. Each word may appear several 
times in the factorization of $p$
and also in the factorization of different aperiodic directed closed walks.
 Let $B(D')$ be the set-system of all the words (with repetition) appearing in the factorizations
of the aperiodic directed closed walks of $D'$. 

It directly follows from Lemma~\ref{lem:coin}, the lemma on coin arrangements, that the sum of all monomial summands $S$ in the expansion of $\prod_{p\in A_1}(1+W(p))$, which have the same set of 'coins' of $B(D')$ with more than 
one element, is zero. Hence the monomial summands $S$ which survive in the expansion of 
$\prod_{p\in A_1}(1+W(p)$  all have the set of coins of $B(D')$ consisting of exactly one word. Such $S$ cannot
have $a_1$ with exponent bigger than one. Now we can repeat the same consideration for the other
edges different from $e(v), v\in V$.

Hence the only terms of the expansion of the infinite product that survive have all $M_e, e\neq e(v)$,
with the exponent at most one. 
 
We know from Observation \ref{o.d} that the collections of the edge-disjoint directed cycles 
of length at least three in $D$ correspond to the collections of the aperiodic directed closed walks of $D'$ where each edge of $D'$ and each vertex $v_g$ appear at most once; by above, these
exactly {\it have chance} to survive in the infinite product. 

Each term of $\xi(M)$ may be expressed several times as a product of aperiodic directed closed walks of $D'$,
\emph{but} only one such expression survives in the infinite product since if a closed walk goes
through an edge $e(v)$ or through a vertex $v_g$ more than once, its rotation is defined to be zero.
Hence $\xi(M)$ is counted correctly in the infinite product. 
  \end{proof}

\smallskip
\begin{remark}
\label{rem.ihz}
\em{
Let us write $\rho(p)= (-1)^{\rot(p)}$. 
Without the zero values of $\rho(p)$, the function $\rot(p)$ is additive when we 'smoothen' $p$
into directed cycles.  The integer lattice generated by the directed cycles has a basis
which may be constructed e.g., from the ear-decomposition~\cite{GL}; the function $\rot(p)$ may
be split into contributions of the edge-transitions for the basis, and since it is a basis, 
it may be split into contributions of the edge-transitions also for all the directed cycles. 
Hence {\bf if} the additivity property holds, $\rot(p)$ may be split into the contributions $\rot (t)$ of the edge-transitions $t$ for the aperiodic closed walks. Hence we would have 
 $$
 \prod_p (1- \rho(p)w(p))= \prod_p (1- \prod_{t\text{ transition of }p}(-1)^{\rot (t)}w(t)).
 $$
This formula transforms the infinite product into the {\it Ihara-Selberg function}. It was studied  by Bass in \cite{Bass} who proved that it is equal to a determinant. A combinatorial proof was given by Foata and Zeilberger in \cite{FZ}. 
 
Having the zero values of $\rho(p)$ it is not clear how to split the rotation to individual edge-transitions. A determinant-type formula may however exist.
}
\end{remark}

\subsection{Acknowledgement}
The authors thank Mireille Bousquet-Melou, Marie-Line Chabanal, Philippe Flajolet, Bojan Mohar, Gilles Schaeffer, and Alexander Zvonkin for helpful discussions. The authors also thank the referees for their comments and suggestions to improve this paper and for the literature on rigorous treatment of the matrix integral method in particular. 

Martin Loebl gratefully acknowledges the support of CONICYT via grant Anillo en Redes.
Mihyun Kang is supported by the Deutsche Forschungsgemeinschaft (DFG  Pr 296).

\end{document}